\documentclass[a4paper,12pt]{article} 

\usepackage[latin1]{inputenc} 
\usepackage[T1]{fontenc} 
\usepackage{indentfirst} 
\usepackage{amsfonts} 
\usepackage{amssymb}
\usepackage{amsmath} 
\usepackage{amsopn} 
\usepackage{amscd} 
\usepackage{amsthm} 
\usepackage{graphicx} 
\usepackage{tikz} 

\usetikzlibrary{matrix} 

\title{A Small Model for the Cohomology of Some Principal Bundles}
\author{Samuel Tinguely}
\date{}

\newtheorem{thm}{Theorem}
\newtheorem{prop}{Proposition}
\newtheorem{cor}{Corollary}
\newtheorem{lem}{Lemma}

\newcommand{\ggoth}[0]{\mathfrak{g}}

\newcommand{\R}{\mathbb{R}}

\newcommand{\APL}{A_{PL}}
\newcommand{\QIto}{\stackrel{\simeq}{\to}}
\newcommand{\QIfrom}{\stackrel{\simeq}{\leftarrow}}

\newcommand{\Def}{{\bfseries Definition :} }
\newcommand{\Pf}{{\bfseries Proof :} }

\newcommand{\Akn}{{\bfseries Aknowledgements :} }

\DeclareMathOperator{\tr}{tr}
\DeclareMathOperator{\Image}{Im}

\DeclareMathOperator{\diff}{d}

\begin{document}

\maketitle

\begin{abstract}
Let $G$ be a compact, connected and simply connected Lie group, and $\Omega G$ the space of the loops in $G$ based at the identity.
This note shows a way to compute the cohomology of the total space of a principal $\Omega G$-bundle over a manifold $M$, from the cohomology of $G$, the differential forms on $M$ and the characteristic classes of the bundle.
The equivariant situation is also treated.
\end{abstract}

\section{Introduction}

Let $P \to M$ be a principal $G$ bundle, where $G$ is a reductive Lie group and $M$ a differential manifold. Let us recall that the real cohomology ring of $G$ is an exterior algebra over a graded vector space $V_G$, with basis $\{x_i\}$ : $H^*(G)\cong \Lambda V_G$ (in this paper, $\Lambda V$ denotes the free graded-commutative algebra over the graded vector space $V$). Now if $BG$ is a classifying space for $G$, then $H^*(BG)\cong \Lambda V_B$, where $V_B=V_G[-1]$ denotes the vector space isomorphic to $V_G$, but with grading shifted by one, so that there is a base $\{y_i\}$ of $V$ such that $\deg y_i = \deg x_i+1$ and a map $\tau: V_G \to V_B; x_i \mapsto y_i$ called the \emph{distinguished transgression}.  Now for any given principal $G$-bundle $P \to M$, any choice of a connection provides a way to pull the generators of $H^*(BG)$ (that is, the elements of $V_B$) back to $\Omega^*(M)$ (that is, the differential forms on $M$). Such a pullback of the image of an element $x_i$ of $V_G$ via the distinguished transgression is a representative of the \emph{characteristic class} corresponding to $x_i$, and denoted by $c_i$ (the class and the representative are denoted the same way).

Greub, Halperin and Vanstone prove the following 
\begin{thm} Let $P\to M$ be a principal $G$-bundle over a smooth manifold $M$. The real cohomology of $P$ is isomorphic to the homology of the complex $\Omega^*(M)\otimes \Lambda V_B$ with the differential defined as follows :
\[\diff(\omega\otimes1) = \diff_{dR}\omega\otimes1\]
and\[
\diff(\omega\otimes x_{1}\wedge...\wedge x_r) = \diff_{dR}\omega\otimes x_1\wedge...\wedge x_r + (-1)^{|\omega|}\sum_{j=1}^r(-1)^j c_j\cdot\omega\otimes x_1\wedge...\widehat{x_j}...\wedge x_r
\]
\qed\end{thm}


Using tools from rational homotopy theory, Félix, Halperin and Thomas proved a very similar result but for a much wider class of groups and base spaces : they only require $G$ to be a path connnected topological group with finite dimentional rational homology, and $M$ to be a simply connected CW-complex (see chapter 15, § (f) in \cite{FHT}).

These tools allow us to turn our interest towards another class of principal bundles.
Now $G$ is a compact, connected and simply connected Lie group. If we denote by $\Omega_e G$ the space of loops in $G$, based at the identity, and the space of paths in $G$ starting at the identity by $P_eG$, every $\Omega_e G$-bundle can be pulled back from the path-loop fibration $\Omega_e G \to P_eG \to G$.
If $M$ is a connected, simply connected manifold, let $f:M\to G$ be the classifiying map for the principal $\Omega_e G$-bundle $\Omega_e G \to P \to M$ and $c_i=f^*\tilde{x}_i$, where $\tilde{x}_i$ denotes the invariant representative of the cohomology class $x_i$. Then the following holds true.
\begin{thm}
Let $P\to M$ be a principal $\Omega_eG$-bundle over a smooth compact, connected and simply-connected manifold $M$. The cohomology of $P$ is given by the homology of the complex $\Omega^*(M)[y_1,... ,y_r]$, where the $y_i$'s are generators corresponding to those of $H^*(G)$, but of degree one less, endowed with the differential
\[
\diff(\omega p(\mathbf{y})) = \diff(\omega)P(\mathbf{y}) + (-1)^{|\omega|} \sum_{i=1}^r \omega \wedge c_i \frac{\partial}{\partial y_i}p(\mathbf{y})
\]
\end{thm}
Note that the choice of $c_i$'s is not unique, as will be made precise in Section \ref{Sec:Result}.

With not so much more effort, one can prove an equivariant version of this. First recall that the equivariant cohomology algebra of $G$ acting on itself by conjugation is exactly $\Lambda (V_G \oplus V_B)$ (see Lemma \ref{HGG}). If $p:P \to M$ is the pull-back of the universal principal $\Omega_e G$-bundle by an equivariant map $f:M\to G$ with respect to the conjugation action, and if we choose closed elements $c_i(\xi)$ of the Cartan model $C_G(M)$ for the equivariant cohomology of $M$ such that $c_i(\xi)$ represents $f_G^*(x_i)$ (the pullback along $f$ in equivariant cohomology), then we have the following
\begin{thm}
The $G$-equivariant cohomology of $P$ is isomorphic, as an algebra, to the cohomology of the complex
$C_G(M)[y_1,... ,y_r]$ with differential
\[
\diff : w(\xi) p(\mathbf{y}) \mapsto \diff_{DR}(w(\xi))p(\mathbf{y}) - \iota_\xi w(\xi)p(\mathbf{y}) + (w\wedge c_i)(\xi) \sum_{i=1}^r \frac{\partial}{\partial y_i}p(\mathbf{y})
\]
\end{thm}
This result is rather interesting since the equivariant cohomology of such a space offers a control of the cohomology of the symplectic quotient $P // LG$ (see \cite{BTW} and \cite{AMM}). E. Meinrenken mentioned this model in his lecture notes on the subject \cite{Mein}.

\Akn This work was supported in part by the grant 140985 of the Swiss Nation Science Foundation. The author would like to thank A. Alekseev, M. Franz, J.-C. Hausmann and K. Hess for many fruitful discussions, as well as J. Huebschmann and E. Meinrenken for their helpful coments.

\section{The tools}\label{tools}

\subsection{Polynomial differential forms}\label{PolyDiff}

In this section, we recall the construction of polynomial differential forms, and one of the main results of interest for us.

We call \emph{cochain algebra} a differential graded algebra concentrated in non-negative degrees, with a differential of degree $+1$.

For a topological space $X$, let $S_*(X)$ be the set of singular simplices in $X$, seen as a simplicial set (see e.g. \cite{May}, example 1.5).

Let us consider the free graded commutative algebra $\Lambda(t_0,...,t_n,y_0,...,y_n)$ (with real coefficients) where the basis elements $t_i$ have degree 0 and $y_j$ have degree 1. There is a unique derivation in this algebra specified by $t_i\mapsto y_i$ and $y_i \mapsto 0$. This derivation preserves the ideal $I_n$ generated by the two elements $\sum_0^nt_i-1$ and $\sum_0^ny_j$, so we can define the quotient differential algebra
\[
(\APL)_n=\frac{\Lambda(t_0,...,t_n,y_0,...,y_n)}{(\sum t_i-1,\sum y_j)}
\]
\[
\diff t_i = y_i \mbox{ and } \diff y_i = 0
\]

The cochain algebra morphisms 
$\partial_i:(\APL)_{n+1}\to(\APL)_n$ and $s_j:(\APL)_{n}\to(\APL)_{n+1}$ 
uniquely specified by
\[
\partial_i : t_k \mapsto \left\{\begin{array}{ll}t_k & , k<i \\ 0 & , k=i \\ t_{k-1} & ,k>i\end{array}\right.
\mbox{ and }
s_j : t_k \mapsto \left\{\begin{array}{ll}t_k & , k<j \\ t_k+t_{k+1} & , k=j \\ t_{k+1} & ,k>j\end{array}\right.
\]
give $\APL = \bigoplus_{n\geq0}(\APL)_n$ a structure of simplicial cochain algebra.

We can now define the functor $\APL$ from the category of topological spaces to the category of cochain algebras :
$\APL^p(X)$ is the set of simplicial set morphism from $S_*(X)$ to $\APL^p$ (the set of elements of degree $p$ in $\APL$), with point-wise addition, scalar and internal multiplication, as well as differentiation. If $f:X\to Y$ is a continuous map, then $\APL(f):\APL(Y)\to\APL(X)$ is the morphism of cochain algebra defined by precomposition by $S_*(f)$.

For any topological space $X$, $\APL(X)$ is \emph{weakly equivalent} to $C^*(X)$, that is there are commutative cochain algebras $(C(0),\diff),...,(C(k),\diff)$ such that 
\[\APL(X) \stackrel{\simeq}{\rightarrow} (C(0),\diff) \stackrel{\simeq}{\leftarrow} \cdots \stackrel{\simeq}{\rightarrow} (C(k),\diff) \stackrel{\simeq}{\leftarrow} C^*(X)\]
If $M$ is a smooth manifold, then $\APL(M)$ is also {weakly equivalent} to $\Omega^*(M)$.

The elements of the cochain algebra $\APL(X)$ are called \emph{polynomial differential forms}.

\subsection{Sullivan algebras and Sullivan models}\label{Sull}

The main tools of this paper are Sullivan algebras, Sullivan models, and their relative counterpart. In this section we expose briefly their definitions and some of their properties, following \cite{FHT}. Unless stated otherwise, all results come from this book, and all algebras have a unit.

\Def A \emph{relative Sullivan algebra} is a commutative cochain algebra of the form $(B\otimes\Lambda V,\diff)$, where
\begin{itemize}
\item $(B,\diff) = (B\otimes 1,\diff)$ is a sub cochain algebra, and $H^0(B)=\R$
\item $1\otimes V = V = \bigoplus_{p\geq1}V^p$ and $\Lambda V$ is the free commutative algebra on $V$
\item $V=\bigcup_{k=0}^\infty V(k)$, where $V(0) \subset V(1) \subset \cdots$ is an increasing sequence of graded subspaces such that
\[\diff:V(0)\to B \mbox{ and } \diff:V(k)\to B\otimes \Lambda(V(k-1)), k\geq 1\]
\end{itemize}

We identify $B=B\otimes1$ and $\Lambda V=1\otimes\Lambda V$. The sub-cochain algebra $(B,\diff)$ is called the \emph{base algebra} of $(B\otimes\Lambda V,\diff)$.


Now let $\varphi:(B,\diff)\to(C,\diff)$ be a morphism of commutative cochain algebras with $H^0(B)=\R$.

\Def A \emph{Sullivan model} for $\varphi$ is a quasi-isomorphism of cochain algebras
\[m:(B\otimes\Lambda V,\diff)\QIto(C,\diff)\]
such that $(B\otimes\Lambda V,\diff)$ is a relative Sullivan algebra with base $(B,\diff)$ and $m|_B=\varphi$.

If $f:X\to Y$ is a continuous map then a Sullivan model for $\APL(f)$ is called a \emph{Sullivan model for $f$}. 


\Def A Sullivan algebra or a Sullivan model is called \emph{minimal} if
\[\Image \diff \subset B^{+}\otimes\Lambda V + B\otimes\Lambda^{\geq2}V\]

The special case where $(B,\diff)=\R$ and $\varphi:(B,\diff)\to(C,\diff)$ is the canonical morphism $\R \to (A,\diff) ; 1\mapsto1$ gets particular attention :

{\bfseries Definitions :} 
\begin{enumerate}
    \item A relative Sullivan algebra with base $\R$ is simply called a \emph{Sullivan algebra}.
    \item A Sullivan model for $\R \to (C,\diff) ; 1\mapsto1$ is called a \emph{Sullivan model for $(C,\diff)$}.
    \item If $X$ is a path-connected topological space, a Sullivan model for $\APL(X)$ is called a \emph{Sullivan model for $X$}.
    \item A Sullivan algebra $(\Lambda V,\diff)$ is called \emph{minimal} if
	   \[\Image \diff \subset \Lambda^{\geq2}V.\]
\end{enumerate}

We can assure the existence of a Sullivan model under some hypotheses :

\begin{prop}\label{RelSullExists}
A morphism $\varphi:(B,\diff)\to(C,\diff)$ of commutative cochain algebras admits a Sullivan model if $H^0(B)=\R=H^0(C)$ and $H^1(\varphi)$ is injective.

Moreover, any commutative cochain algebra $(A,\diff)$ with $H^0(A)=\R$ (and any path-connected topological space) admits a \emph{unique} minimal Sullivan model.
\qed
\end{prop}

We can combine Propositions 12.8 and 12.9 in \cite{FHT} to get the following

\begin{prop}\label{Lifting}
Let $\eta : (A,\diff)\QIto(C,\diff)$ be a quasi-isomorphism of commutative cochain algebras. Let $(\Lambda V,\diff)$ be a Sullivan algebra, and $\psi : (\Lambda V,\diff)\to(C,\diff)$ a morphism of cochain algebras. Then there is a morphism of commutative cochain algebras $\varphi : (\Lambda V,\diff)\to(A,\diff)$ such that $H(\eta\circ\varphi)=H(\psi)$. Any two such morphisms  $\varphi_1, \varphi_2 : (\Lambda V,\diff)\to(A,\diff)$ satisfy $H(\varphi_1) = H(\varphi_2)$
\qed
\end{prop}

\begin{equation}\label{LiftingDiag}
\begin{tikzpicture}
\matrix [matrix of math nodes,row sep=1cm, column sep=1.5cm]
{
 & |(A)|(A, \diff) \\
|(V)|(\Lambda V, \diff) & |(C)|(C, \diff) \\
};
\tikzstyle{every node}=[midway,font=\small]
\draw[-stealth] (V) -- (C) node [above] {$\psi$};
\draw[dashed, -stealth] (V) -- (A) node [above] {$\varphi$};
\draw[-stealth] (A) -- (C) node [right] {$\eta$} node  [left] {$\simeq$};
\end{tikzpicture}
\end{equation}

It is called the \emph{lifting lemma} because it say that for any given $\psi$ and $\eta$ as in Diagram \ref{LiftingDiag}, there exists a $\varphi$ making the diagram commute.

As it stands, this proposition isn't so useful to us, but it has two crucial corollaries :

\begin{cor}\label{SullWeakEq}
Let $(A,\diff)$ and $(B,\diff)$ be two commutative cochain algebras, and $\psi : (\Lambda V,\diff)\QIto(A,\diff)$ a Sullivan model for $(A,\diff)$. If $(A,\diff)$ and $(B,\diff)$ are weakly equivalent, then there exists a Sullivan model $\varphi : (\Lambda V,\diff)\QIto(B,\diff)$ for $(B,\diff)$.
\end{cor}
\Pf Recall that for $(A,\diff)$ and $(B,\diff)$ to be weakly equivalent means that there are cochain algebras $(C(0),\diff),...,(C(k),\diff)$ so that 
\[(A,\diff) \stackrel{\simeq}{\rightarrow} (C(0),\diff) \stackrel{\simeq}{\leftarrow} \cdots \stackrel{\simeq}{\rightarrow} (C(k),\diff) \stackrel{\simeq}{\leftarrow} (B,\diff).\]
Now use Diagram \ref{LiftingDiag} to conclude.
\qed

\begin{cor}\label{SullRepr}
Let $(A,\diff)$ and $(B,\diff)$ be two commutative cochain algebras, and let $m_A:(\Lambda V_A,\diff)\to(A,\diff)$ and $m_B:(\Lambda V_B,\diff)\to(B,\diff)$ be Sullivan models for $A$ and $B$ respectively. Then for any morphism $f:(B,\diff)\to(A,\diff)$ there is a \emph{Sullivan representative} for $f$, that is a morphism $\phi:\Lambda V_B \to \Lambda V_A$ such that $H(f \circ m_B) = H(m_A \circ \phi)$.
\qed\end{cor}

Now suppose $(B\otimes\Lambda V,\diff)$ is a relative Sullivan Algebra and $\psi:(B,\diff)\to(B',\diff)$ is a morphism of commutative cochain algebras with $H^0(B')=\R$. Then $\psi$ gives $(B',\diff)$ a structure of $(B,\diff)$ module, and the cochain algebra
\[(B',\diff)\otimes_{(B,\diff)}(B\otimes\Lambda V,\diff)=(B'\otimes\Lambda V,\diff)\]
is a relative Sullivan algebra with base $(B',\diff)$. It is called the \emph{pushout of $(B\otimes\Lambda V,\diff)$ along $\psi$}. 

This pushout construction has the really useful property that it preserves quasi-isomorphisms :

\begin{lem}\label{QIso}
If $\psi$ is a quasi-isomorphism, then so is $\psi\otimes id : (B\otimes\Lambda V,\diff) \to (B'\otimes\Lambda V,\diff)$.
\qed
\end{lem}

If $\varepsilon:B\to\R$ is any augmentation, pushing $(B\otimes\Lambda V,\diff)$ out along $\varepsilon$ yields a Sullivan algebra $(\Lambda V,\bar{\diff})$, which is called the \emph{Sullivan fibre at $\varepsilon$}.

\subsection{Application to fibrations}\label{fibrations}

Let $p:X\to Y$ be a fibration with fibre $F$, and $q:Z\to A$ be the pullback of $p$ along a continuous mapping $f:A\to Y$, where both $A$ and $Y$ are simply connected. 
\begin{equation}
\begin{CD}
Z @>\bar{f}>> X \\
@VqVV @VpVV \\
A @>f>> Y
\end{CD}
\end{equation}
Assume in addition that one of $H_*(X)$, $H_*(Y)$, $H_*(Z)$, $H_*(A)$ or $H_*(F)$ is of finite type.
Let us choose Sullivan models $m_Y:(\Lambda V_Y,\diff)\to\APL(Y)$, $n_A:(\Lambda W_A,\diff)\to\APL(A)$ and $m:(\Lambda V_Y\otimes\Lambda V,\diff)\to\APL(X)$ (modelling $p\circ m_Y$), as well as a Sullivan representative for $f$
\[\psi:(\Lambda V_Y,\diff)\to(\Lambda W_A,\diff)\]

Then one can define the morphism
\[
\xi:(\Lambda W_A,\diff)\otimes_{(\Lambda V_Y,\diff)}(\Lambda V_Y\otimes\Lambda V,\diff)\to\APL(Z)
\]
to be the composition
\begin{multline}
(\Lambda W_A,\diff)\otimes_{(\Lambda V_Y,\diff)}(\Lambda V_Y\otimes\Lambda V,\diff)
\stackrel{n_A \otimes m}{\longrightarrow}
\APL(A)\otimes\APL(X) \\
\stackrel{\APL(q) \otimes \APL(\bar{f})}{\longrightarrow}
\APL(Z)\otimes\APL(Z)
\stackrel{\cdot}{\to}\APL(Z)
\end{multline}
where $\cdot$ is simply the multiplication in $\APL(Z)$.
We rewrite $(\Lambda W_A,\diff)\otimes_{(\Lambda V_Y,\diff)}(\Lambda V_Y\otimes\Lambda V,\diff)$ as $(\Lambda W_A\otimes\Lambda V,\diff)$. One can show the following

\begin{prop}\label{15.8}
$\xi:(\Lambda W_A\otimes\Lambda V,\diff)\to\APL(Z)$ is a Sullivan model for $Z$.
\qed
\end{prop}

Morally, this means that \emph{the pushout of the models is a model for the pullback}.

Note that if $A$ is a point $y_0$ of $Y$ and $f$ is the inclusion (so that $Z=F$), then the minimal Sullivan model for $A=y_0$ is $(\R,0)\QIto\APL(A)$, so that $(\Lambda W_A\otimes\Lambda V,\diff)=(\Lambda V,\bar{\diff})$, where $\bar{\diff}$ comes from the pushout along the aumentation $\psi:(\Lambda V_Y,\diff)\to\R$ induced by $f$. We can now formulate the following 
\begin{cor}\label{fibre}
If $Y$ is simply connected and $X\to Y$ is a fibration with fibre $F$, then
\[\xi : (\Lambda V,\bar{\diff})\QIto \APL(F)\]
is a Sullivan model for $F$.
\qed\end{cor}
which can be rephrased as `\emph{the fibre of the model is a model for the fibre}'.

\section{The result}\label{Sec:Result}

Let $G$ be a simply connected Lie group, and $P_eG \stackrel{\pi}{\rightarrow} G$ the path fibration over $G$. In this situation, we have a map of commutative cochain algebras $\APL(\pi):\APL(G) \to \APL(P_eG)$.

The cohomology ring of $G$ is the free graded commutative algebra $\Lambda V_G$, where $V_G$ is a finite dimentional graded vector space concentrated in odd degrees, with $V_G^1=0$. Moreover, $H^*(G) = (\Lambda V_G,0)$ is a Sullivan algebra, so $H^*(G)\stackrel{\cong}{\to}\Omega^*(G)^G\stackrel{\simeq}{\hookrightarrow}\Omega^*(G)$ is a Sullivan model, and since $\APL(G)$ and $\Omega^*(G)$ are weakly equivalent, their Sullivan models are identified, so there is a Sullivan model $m_G:(\Lambda V_G,0)\QIto\APL(G)$. Now let $\{x_i\}_{i=1}^r$ be a basis for $V_G$, and $V$ be a vector space isomorphic to $V_G$ with grading shifted by one, and basis $\{y_i\}_{i=1}^r$ such that $\deg y_i = \deg x_i-1$. Let us call this correspondance $\delta : V \to V_G$ and extend it to a derivation $\diff$ on $\Lambda(V_G\oplus V)$. Since $\delta$ is an isomorphism of vector spaces and $V^0=0$, 
$(\Lambda(V_G\oplus V),\diff)$ is acyclic. Now let us push this out along $m_G$ to get $(\APL(G)\otimes \Lambda V,\diff)$. Since $m_G$ is a quasi-isomorphism, $m_G\otimes id : (\Lambda(V_G\oplus V),\diff) \to (\APL(G)\otimes \Lambda V,\diff)$ is also a quasi-isomorphism, by Lemma \ref{QIso}.

Let us now define $m:(\APL(G)\otimes\Lambda V,\diff) \to (\APL(P_eG),\diff)$ such that $m$ coincides with $\APL(\pi)$ on $\APL(G)$. Since $P_eG$ is contractible, every cocycle of degree $\geq1$ is a coboundary. In particular, $m\circ m_G(x_i)=\diff a_i$ for some $a_i\in\APL(P_eG)$. We now define $m(y_i)$ to be $a_i$. These assignments extend to a morphism of cochain algebra $m:(\APL(G)\otimes\Lambda V,\diff) \to (\APL(P_eG),\diff)$ that that restricts to $\APL(\pi)$ on $\APL(G)$.

\begin{lem}
$m:(\APL(G)\otimes\Lambda V,\diff) \to (\APL(P_eG),\diff)$ is a Sullivan model for $\APL(\pi):\APL(G) \to \APL(P_eG)$.
\end{lem}
\Pf
By construction, $m$ is a morphism of differential graded commutative algebras extending $\APL(\pi)$. Since both $\APL(P_eG)$ and $(\APL(G)\otimes\Lambda V,\diff)$ are acyclic, $m$ is a quasi-isomoprhism. Now we only need to show that $(\APL(G)\otimes\Lambda V,\diff)$ is a relative Sullivan algebra.

By construction, $(\APL(G),\diff)$ is a sub cochain algebra and $V=\oplus_{p\geq1}V^p$, and by connectedness of $G$, $H^0(\APL(G))=\R$. We only have to check the nilpotence condition, which is trivially satisfied by setting $V(0)=V$, since $\diff$ maps $V$ to $\APL(G)$.
\qed

An immediate consequence of this and Corollary \ref{fibre} is that there is a quasi-isomorphism $(\Lambda V, 0)\to \APL(\Omega_e G)$, which means that we can formulate the following

\begin{cor}
If $G$ is a connected, simply-connected Lie group, and if $\{x_i\}$ are the generators of its cohomology algebra, then the cohomology algebra of $\Omega_e G$ is the polynomial ring $\R[\{y_i\}]$, where $\deg(y_i)=\deg(x_i)-1$.
\qed
\end{cor}

Let us now examine our situation :
\begin{equation}\label{diagbase}
\begin{CD}
P   @>\bar{f}>>      P_eG   \\
@V\pi_PVV     @V\pi VV  \\
M   @>f>>     G
\end{CD}
\end{equation}
The quasi-isomorphism $m_G:(\Lambda V_G,0) \to \APL(G)$ is the minimal model for $G$. Let us choose a minimal Sullivan model $n_M:(\Lambda W_M, \diff) \to \APL(M)$. Then by Proposition \ref{15.8}, if $H^*(\pi_P)$ is injective (e.g. if $M$ is simply connected), we have the following

\begin{prop}
The morphism
\[
\xi = \APL(\pi_P)n_M\cdot\APL(\bar{f})m : (\Lambda W_M,\diff)\otimes_{\Lambda V_G}(\Lambda V_G \otimes \Lambda V, \diff) \to \APL(P)
\]
is a Sullivan model for $P$
\end{prop}
We will most often write $(\Lambda W_M \otimes \Lambda V, \diff)$ for $(\Lambda W_M,\diff)\otimes_{\Lambda V_G}(\Lambda V_G \otimes \Lambda V, \diff)$. 

Now since $\APL(M)$ and $\Omega^*(M)$ are weakly equivalent, their Sullivan models are the same, hence there is a quasi-isomorphism $(\Lambda W_M,\diff) \QIto \Omega^*(M)$. Using Lemma \ref{QIso} one can replace $(\Lambda W_M,\diff)$ by  $\Omega^*(M)$ in $\xi$ to get a weak equivalence between $\APL(P)$ and $(\Omega^*(M) \otimes \Lambda V, \diff)$. The differential on this complex is the following : if $n_{DR}:(\Lambda W_M,\diff)\QIto\Omega^*(M)$ is the Sullivan model for $\Omega^*(M)$ corresponding to $n_M:(\Lambda W_M, \diff) \to \APL(M)$, $\psi:(\Lambda W_M, \diff)\to(\Lambda W_M,\diff)$ is a Sullivan representative for $f$, and we write $c_i =  n_{DR} \psi (x_i)$, we have that :
\[
\diff : y_i = 1 \otimes y_i \mapsto c_i \ \ \ \ \   \diff : w = w \otimes 1 \mapsto \diff w \otimes 1
\]
for any element $w$ of $\Omega^*(M)$. 
We thus have proved the following 
\begin{thm}\label{Result}
In the situation of the diagram \ref{diagbase}, the cohomology of $P$ is isomorphic, as an algebra, to the cohomology of the complex
\[
\Omega^*(M)[y_1,...y_r]
\]
with differential
\[
\diff : w P(\mathbf{y}) \mapsto \diff w P(\mathbf{y}) + \sum_{i=1}^r w\wedge c_i \frac{\partial}{\partial y_i}P(\mathbf{y})
\]
where the $c_i$'s are defined as above, and $\deg y_i = \deg x_i-1$.
\qed
\end{thm}

Note that the resulting graded algebra depends neither on the choice of the Sullivan representative for $f$, nor of the choice of the model $n_{DR}:(\Lambda W_M,\diff)\QIto\Omega^*(M)$;
the only important thing is that $n_{DR}\circ\psi(x_i)$ represents the same cohomology class as $\APL(f)\circ m_G(x_i)$.
Therefore, one can choose for instance to map $x_i$ to the pullback by $f$ of the bi-invariant representatives of $x_i$.

Of course, when $M$ is formal, we can go one step further : since $H^*(M)$ is weakly equivalent to $\APL(M)$, a Sullivan model for $\APL(M)$ is also one for $H^*(M)$, so we can pushout $(\Lambda W_M \otimes \Lambda V,\diff)$ to $(H^*(M) \otimes \Lambda V,\diff)$, which is still a model for $\APL(P)$. In other words :
\begin{thm}\label{ResultFormal}
In the situation of the diagram \ref{diagbase}, and if $M$ is formal, the cohomology of $P$ is isomorphic, as an algebra, to the cohomology of the complex
\[
H^*(M)[y_1,...y_r]
\]
with differential
\[
\diff : w P(\mathbf{y}) \mapsto \sum_{i=1}^r w\wedge c_i \frac{\partial}{\partial y_i}P(\mathbf{y})
\]
where the $c_i$'s are the characteristic classes of $P$, and $\deg y_i = \deg x_i-1$.
\qed
\end{thm}

\section{An example}

These results allow for computations of some not entirely trivial bundles. For instance, let $G$ be $SU(3)$. In this case, $H^*(SU(3))=\bigwedge[\eta_3,\eta_5]$, with $\deg \eta_3 = 3$ and $\deg \eta_5 = 5$. 
Let us then define
\[M = \{g\bar{g}^{-1} | g \in SU(3)\} = \{gg^\top | g \in SU(3)\}.\]

Let us note that $M$ is the orbit of the identity by the 'twisted' action $h \mapsto gh\bar{g}^{-1}$, and a quick calculation shows that the stabilizer of the identity for this action is $SO(3)$, so that $M=SU(3)/SO(3)$. Of course, the restriction of the biinvariant representative of $\eta_3$ in $\Omega^*(SU(3))$ to $SO(3)$ is still biinvariant. To see this, write $\eta_3$ in matrix form $\tr (g^{-1}\diff g)^3$ and see $SO(3)$ as a subspace of $SU(3)$ (see e.g. \cite{GHV}, chap. 6, § 7. for a justification). This implies that $i^*:H^*(SU(3)) \to H^*(SO(3))$ is surjective, so one can deduce (see e.g. \cite{GHV} section 11.6) that $H^*(M) = \R[x_5]/(x_5^2)$.

In particular, $M$ is formal : any choice of a representative of the class $x_5$ yields a quasi-isomorphism $H^*(M) \QIto \Omega^*(M)$.

Now in order to be able to apply our result, we have to know $f^*$. Trivially, $f^*\eta_3 =0$, so we only have to examine $f^*\eta_5$. Now the map $SU(3) \to SU(3)$, $g \mapsto g\bar{g}^{-1}$ is non-zero in degree $5$ cohomology, so $f^*\eta_5$ cannot be zero, which means that $f^*\eta_5 = x_5$.

Let now $P$ be the pull back of the path fibration of $SU(3)$ by the inclusion $f:M \to SU(3)$. Since $M$ is formal, we can use Theorem \ref{ResultFormal} :
\[
H^*(P) = H^*\left(H^*(M)[y_2, y_4], \diff = f^*x_3\frac{\partial}{\partial y_2}+f^*x_5\frac{\partial}{\partial y_4}\right)
\]
Where $\deg y_k = k$.

Let us now do concrete computations. The cocycles are linear combinations of expressions of the form $y_2^p$ or $\eta_5\otimes y_2^p y_4^q$. The only coboudaries are the $\diff(y_2^p y_4^q) = \eta_5\otimes y_2^p y_4^{q-1}$, which includes all $\eta_5\otimes y_2^p y_4^q$. The cohomology of $P$ is therefore simply $\R[y_2]$.

\section{The equivariant case}

Of particular interest are the cases where $G$ acts on $M$, and the map $p:M\to G$ is equivariant with respect to the adjoint action of $G$ on itself. The action of $G$ by conjugation on $P_eG$ descend to conjugation on $G$ via the projection to the endpoint $\pi:P_eG \to G$, so that $\pi$ is a universal $G$-equivariant $\Omega_eG$-bundle.


This implies that an equivariant map $M\to G$ induces a pullback diagram, that is entirely composed of $G$-equivariant mappings :
\begin{equation}
\begin{CD}
P @>>> P_eG\\
@VVV @VVV\\
M @>>> G
\end{CD}
\end{equation}
Since the Borel construction is functorial, we can compute the $G$-equivariant cohomology from topological spaces as in this diagram :
\begin{equation}
\begin{CD}
P_G @>>> P_eG_G\\
@VVV @VVV\\
M_G @>>> G_G
\end{CD}
\end{equation}
One checks easily that this is still a pullback diagram, so that we can use the results of Section \ref{fibrations}.

Since we still have actual topological spaces at hand, we still have that the pushout of the models are a model for the pullback. So first, we have to establish a model for $\APL(G_G)$ :

\begin{lem}\label{HGG}
$m_{G_G}:(\Lambda (V_{G}\oplus V),0)\QIto\APL(G_G)$ is a Sullivan model for $G_G$.
\end{lem}
\Pf First remark that the bundle $G \to G_G \to BG$ has a section (because the identity of $G$ is a fixed point for the conjugation action), so that $H^*(BG)=H^*_G(pt) \to H^*_G(G) = H^*(G_G)$ is injective. On the other hand, by Proposition \ref{15.8} and the structure of $H^*(G)$ and $H^*(BG)$, there is a Sullivan model of the form $(\Lambda V_{G}\otimes \Lambda V,\diff)\QIto\APL(G_G)$, but since the differentials inside $V_G$ and $V$ are zero, and $H^*(BG) \to H^*(G_G)$ is injective, the global differential in $(\Lambda V_{G}\otimes \Lambda V)$ must be zero : if not, some elements in $\Lambda V = H^*(BG)$ would not be cocylcles, so $H^*(BG) \to H^*(G_G)$ could not be injective.
\qed

Then a construction similar to the one above yields that $(\Lambda (V_{G}\oplus V) \otimes \Lambda V,\diff)$ with essentially the same differential as in the non-equivariant case is a Sullivan model for $\APL(G_G)\to\APL(P_eG_G)$. Now the pushout 
\[
(\Lambda W_{M_G},\diff)\otimes_{(\Lambda (V_{G}\oplus V),0)}(\Lambda (V_{G}\oplus V) \otimes \Lambda V,\diff)
\]
is a model for $H^*(P_G)=H_G^*(P)$, which can be written more concisely as
\[
H_G^*(P)=H^*(\Lambda W_{M_G} \otimes \Lambda V,\diff)
\]

Since $M$ is a manifold, we can compute its equivariant cohomology via the Cartan model (see \cite{GS}, \cite{Car1} and \cite{Car2}). So ideally, we would like to include this model in our calculations. Fortunately, this is possible, because $C^*(M_G)$ is weakly equivlaent to $\APL(M_G)$, and the following

\begin{lem}\label{CartanC*}
The Cartan model is weakly equivalent to $C^*(M_G)$.
\end{lem}
\Pf

This is an easy corollary of a number of results in \cite{GS}. I will simply present the chain of quasi-isomorphism, and indicate the result in this book saying that it is a quasi-isomorphism.

The first few are present in their proof of the equivariant De Rham Theorem, section 2.5. 
$\mathcal{E}$ denotes the set of orthonormal $n$-tuples in $\mathbb{C}^\infty$ (for all values of $n$), which is the inductive limit of $\mathcal{E}_k$, the set of orthonormal $n$-tuples in $\mathbb{C}^\infty$ (for all $n < k$), and $\Omega^*(\mathcal{E})$ denotes the projective limit of the $\Omega^*(\mathcal{E})_k$.
\[C^*((M\times\mathcal{E})/G) \QIfrom \Omega^*((M\times\mathcal{E})/G)\]
\[\Omega^*((M\times\mathcal{E})/G) \stackrel{\cong}{\leftarrow} \Omega^*(M\times\mathcal{E})_{bas}\]
\[\Omega^*(M\times\mathcal{E})_{bas} \QIfrom (\Omega^*(M)\otimes\Omega^*(\mathcal{E}))_{bas}\]
All of these results use the compatibility of cohomology with colimits. The fist line comes from the usual De Rham theorem, the second is well known, and the third is a bit more involved, making use of a spectral sequence argument.

Then we apply Theorem 4.3.1 and its proof, as suggested in section 4.4 (here $W$ denotes the Weyl algebra of $G$) to get :
\[(\Omega^*(M)\otimes\Omega^*(\mathcal{E}))_{bas}  \stackrel{\simeq}{\hookrightarrow} {((\Omega^*(M)\otimes\Omega^*(\mathcal{E}))_{hor})\otimes W)}^G\]
\[{((\Omega^*(M)\otimes\Omega^*(\mathcal{E}))_{hor})\otimes W)}^G \QIfrom (\Omega^*(M)\otimes\Omega^*(\mathcal{E})\otimes W)_{bas} \]
\[(\Omega^*(M)\otimes\Omega^*(\mathcal{E})\otimes W)_{bas} \QIto {((\Omega^*(M)\otimes W)_{hor}\otimes\Omega^*(\mathcal{E}))}^G\]
\[{((\Omega^*(M)\otimes W)_{hor}\otimes\Omega^*(\mathcal{E}))}^G \QIfrom (\Omega^*(M)\otimes W)_{bas}\]
The main ingredients here are the Mathai-Quillen isomorphism, and a filtration argument on the acyclic component. The two central lines are just interchanging the roles of $W$ and $\Omega^*(\mathcal{E})$.

and finally, in section 4.2 we find the last quasi-isomorphism we need :
\[(\Omega^*(M)\otimes W)_{bas} \QIto {(\Omega^*(M)\otimes S(\ggoth^*))}^G\]
Which again comes from the Mathai-Quillen isomorphism, and the fact that $W_{hor} \cong S(\ggoth^*)$.
\qed

So the Sullivan models for $\APL(M_G)$ and ${(\Omega^*(M)\otimes S(\ggoth^*))}^G$ are identified, which implies that there is a quasi-isomorphism $(\Lambda W_{M_G},\diff) \QIto {(\Omega^*(M)\otimes S(\ggoth^*))}^G$ along which we can pushout our model to give
\[
H_G^*(P)=({(\Omega^*(M)\otimes S(\ggoth^*))}^G \otimes \Lambda V,\diff)
\]
and there remains to identify the differential.


Let us write $c_i$ for the image of $x_i$ under the Sullivan representative for $f_G$, followed by the model $(\Lambda W_{M_G},\diff) \QIto {(\Omega^*(M)\otimes S(\ggoth^*))}^G$, and let $y_i$ denote the basis elements of $V$. We can then write :

\begin{thm}
In the situation of the diagram \ref{diagbase}, the $G$-equivariant cohomology of $P$ is isomorphic, as an algebra, to the cohomology of the complex
\[
{(\Omega^*(M)\otimes S(\ggoth^*))}^G[y_1,...y_r]
\]
with differential
\[
\diff : w(\xi) P(\mathbf{y}) \mapsto \diff_{DR}(w(\xi)) - \iota_\xi w(\xi) + \sum_{i=1}^r (w\wedge c_i)(\xi) \frac{\partial}{\partial y_i}P(\mathbf{y})
\]
\qed
\end{thm}


\end{document}